\documentclass{jaca}  

\usepackage{graphicx} 
\usepackage{eqnarray}
\usepackage[all]{xy}
\usepackage{amsmath,amssymb,amscd,amsthm,dsfont,eufrak,euscript,fancyhdr,graphicx,hyperref}
\usepackage{graphics,latexsym,xcolor}

\newtheorem{theo}{Theorem}

\newtheorem{prrop}{Proposition}

\newtheorem{defi}{Definition}

\newcommand{\Rr}{\mathbb{R}}
\newcommand{\Nn}{\mathbb{N}}

\newcommand{\fonction}[5]{\begin{array}[t]{lrcl}#1 :&#2 &\longrightarrow &#3\\&#4& \longmapsto &#5 \end{array}}

\newcommand{\dis}{\displaystyle}

\begin{document}


\title{Variational integrators of fractional Lagrangian systems in the framework of discrete embeddings}



\author{Lo\"ic Bourdin}




\begin{abstract}
This paper is a summary of the theory of discrete embeddings introduced in \cite{bour}. A discrete embedding is an algebraic procedure associating a numerical scheme to a given ordinary differential equation. Lagrangian systems possess a variational structure called Lagrangian structure. We are specially interested in the conservation at the discrete level of this Lagrangian structure by discrete embeddings. We then replace in this framework the variational integrators developed in \cite[Chapter VI.6]{lubi} and in \cite{mars}. Finally, we extend the notion of discrete embeddings and variational integrators to fractional Lagrangian systems.
\end{abstract}

\KeysAndCodes{Lagrangian systems, Variational integrator, Fractional calculus}{70H03, 37K05, 26A33}


\section*{Introduction}
The theoretical framework of embeddings of dynamical systems is initiated by Cresson and Darses in \cite{cres}. A review of the subject is given in \cite{cres3}. An embedding of an ordinary or partial differential equation is a way to give a sense to this equation over a larger set of solutions. As an example, the stochastic embedding developed in \cite{cres} allows to give a meaning of a differential equation over the set of stochastic processes.

We are specially interested in Lagrangian systems covering a large set of dynamical behaviors and widely used in classical mechanics, \cite{arno}. These systems possess a variational structure called Lagrangian structure, i.e. their solutions correspond to critical points of Lagrangian functionals, \cite[p.57]{arno}. The Lagrangian structure is intrinsic and induces strong constraints on the qualitative behavior of the solutions. The conservation of this structure by embedding seems then important. In \cite{cres}, the authors construct stochastic embeddings which preserve the variational structure of Lagrangian systems, i.e. the generalized solutions are also characterized as critical points of generalized Lagrangian functionals. \\

This paper is a summary of the theory of discrete embeddings introduced in \cite{bour} where, as in \cite{cres}, we are interested in the conservation of the Lagrangian structure of Lagrangian systems. We then refer to \cite{bour} for more details and for the proof of some results.

A discrete embedding is an algebraic procedure associating a numerical scheme to a given differential equation, in particular to a given Lagrangian system. On the other hand, defining a discrete embedding induces a discretization of the Lagrangian functional associated and we can develop a discrete calculus of variations on this one: this leads to a numerical scheme called variational integrator. The variational integrators, developed in \cite[Chapter VI.6]{lubi} and \cite{mars}, are then numerical schemes for Lagrangian systems preserving their variational structures.

Thus, we propose the following definition: a discrete embedding is said to be coherent if the two discrete versions obtained (the direct one and the variational integrator) of a Lagrangian system coincide. Hence, a coherent discrete embedding conserves at the discrete level the Lagrangian structure of a Lagrangian system. \\

Recently, many studies have been devoted to fractional Lagrangian systems, \cite{agra}, \cite{cres}. They arise for example in fractional optimal control theory, \cite{fred}. They are difficult to solve explicitly, it is then interesting to develop efficient numerical schemes to such systems. 

Some preliminary results on fractional discrete operators and on the discretization of fractional Euler-Lagrange equations have been discussed by several authors, \cite{agra3}, \cite{bast}, \cite{dubo}. In this paper, we extend the discrete embedding point of view, the corresponding problem of coherence and the associated notion of variational integrator to the fractional case. \\
%

The paper is organized as follows. In Section \ref{section1}, we define the notion of discrete embeddings of differential equations. Section \ref{section2} recalls definitions and results concerning Lagrangian systems and we apply the previous theory of discrete embeddings to Lagrangian systems. Then, we recall the strategy of variational integrators of Lagrangian systems in the framework of discrete embeddings and we finally present the problem of coherence of a discrete embedding. Section \ref{section3} is devoted to the extension of discrete embeddings to the fractional case.

\section{Notion of discrete embeddings}\label{section1}
In this paper, we consider classical and fractional differential systems in $\Rr^d$ where $d \in \Nn ^*$ is the dimension. The trajectories of these systems are curves $q$ in $\mathcal{C}^0 ([a,b],\Rr^d)$ where $a<b$ are two reals. For smooth enough functions $q$, we denote $\dot{q} = dq/dt$ and $\ddot{q} = d^2 q/dt^2$. 

\subsection{Discrete embeddings}\label{section11}

\begin{defi}\label{defde}
Defining a discrete embedding means giving a discrete version of the following elements: the curves $q \in \mathcal{C}^0 ([a,b],\Rr^d)$, the derivative operator $d/dt$ and the functionals $a : \mathcal{C}^0 ([a,b],\Rr^d) \longrightarrow \Rr$. More precisely, it means giving:
\begin{itemize}
\item an application $q \longmapsto q^h $ where $q^h \in (\Rr^d)^{m_1}$, 
\item a discrete operator $\Delta : (\Rr^d)^{m_1} \longrightarrow (\Rr^d)^{m_2}$ discretizing the differential operator $d/dt$,
\item an application $a \longmapsto a^h $ where $a^h : (\Rr^d)^{m_1}  \longrightarrow \Rr $,
\end{itemize}
where $m_1, m_2 \in \Nn^*$.
\end{defi} 

\noindent In order to illustrate Definition \ref{defde}, we define \textit{backward and forward finite differences embeddings}. For all the rest of the paper, we fix $\sigma = \pm$ and $N \in \Nn^*$. We denote by $h = (b-a)/N$ the step size of the discretization and $\tau = (t_k)_{k=0,...,N}$ the following partition of $[a,b]$:
$$ \forall k=0,...,N, \; \; t_{k} = a + k h .$$

\begin{defi}[case $\sigma = -$]
We call backward finite differences embedding denoted by $FDE-$ the definition of the following elements: the application
$$ \fonction{\mbox{disc}}{\mathcal{C}^{0}([a,b], \Rr^d )}{\left( \Rr^d \right) ^{N+1}}{q}{(q(t_k))_{k=0,...,N}} ,$$
and the discrete operator
$$ \fonction{\Delta _-}{\left( \Rr^d \right) ^{N+1}}{\left( \Rr^d \right) ^{N}}{Q = (Q_k)_{k=0,...,N}}{\left( \frac{Q_{k}-Q_{k-1}}{h} \right) _{k=1,..,N}.} $$
\end{defi}

\begin{defi}[case $\sigma = +$]
We call forward finite differences embedding denoted by $FDE+$ the definition of the following elements: the application $\mbox{disc}$ and the discrete operator
$$ \fonction{\Delta _+}{\left( \Rr^d \right) ^{N+1}}{\left( \Rr^d \right) ^{N}}{Q = (Q_k)_{k=0,...,N}}{\left( \frac{Q_{k}-Q_{k+1}}{h} \right) _{k=0,..,N-1}.} $$
\end{defi}

\noindent Let us notice that the discrete analogous of $d/dt$ in $FDE\sigma$ is then $-\sigma \Delta_{\sigma}$. We use these notations in order to be uniform with the fractional notations (see Section \ref{section3}). 

\subsection{Direct discrete embeddings}\label{section12}
Defining a discrete embedding allows us to define a direct discrete version of a given differential equation: 

\begin{defi}
Let be fixed a discrete embedding as defined in Definition \ref{defde} and let \eqref{edo} be an ordinary differential equation of unknown $q \in \mathcal{C}^0([a,b],\Rr^d)$ given by:
\begin{equation}\label{edo}\tag{$E$}
O(q) = 0 ,
\end{equation}
where $O$ is a differential operator shaped as $O = \sum_{i} f_i (\cdot) (d/dt)^i \circ g_i (\cdot) $
where $f_i$, $g_i$ are functions. Then, the direct discrete embedding of \eqref{edo} is \eqref{edoh} the system of equations of unknown $q^h \in (\Rr^d)^{m_1}$ given by:
\begin{equation}\label{edoh}\tag{$E_h$}
O^h(q^h) = 0 ,
\end{equation}
where $O^h$ is the discretized operator of $O$ given by $O^h = \sum_{i} f_i (\cdot)  \Delta^i  \circ  g_i (\cdot) $.
\end{defi}

\noindent As an example, we consider the Newton's equation with friction of unknown $q \in \mathcal{C}^0([a,b],\Rr^d)$ given by:
\begin{equation}\label{ne}\tag{$NE$}
\forall t \in [a,b], \; \; \ddot{q}(t)+\dot{q}(t)+q(t) = 0.
\end{equation}
Then, the direct discrete embedding of \eqref{ne} with respect to $FDE-$ is \eqref{ned} the system of equations of unknown $Q \in (\Rr^d)^{N+1}$ given by:
\begin{equation}\label{ned}\tag{$NE_h$}
\forall k = 2,...,N , \; \; \dis \frac{Q_{k}-2Q_{k-1}+Q_{k-2}}{h^2}+\dis \frac{Q_{k}-Q_{k-1}}{h} +Q_k = 0.
\end{equation}

\noindent The direct discrete embedding of an ordinary differential equation is strongly dependent on the form of the differential operator $O$ (and not on its equivalence class). The process $ O \longrightarrow O^h $ is not an application. For example, the discretized operator $O^h$ of $O = d/dt \circ sin(\cdot) = d/dt (\cdot) \; cos(\cdot)$ is different depending on the writing of $O$. 

\subsection{Direct discrete embeddings of Lagrangian systems}\label{section13}
We recall now classical definitions and theorems concerning Lagrangian systems. We refer to \cite{arno} for a detailed study and for a detailed proof of Theorem \ref{vp}.

\begin{defi}
A Lagrangian functional is an application defined by:
\begin{equation*}
\fonction{\mathcal{L}}{\mathcal{C}^{2}([a,b],\Rr^d)}{\Rr}{q}{\dis \int_{a}^{b} L(q(t),\dot{q}(t),t)dt} 
\end{equation*}
where $L$ is a Lagrangian i.e. a $\mathcal{C}^{2}$ application defined by:
$$ \fonction{L}{\Rr^d \times \Rr^d \times [a,b]}{\Rr}{(x,v,t)}{L(x,v,t).} $$
\end{defi} 
\noindent An \textit{extremal} (or \textit{critical point}) of a Lagrangian functional $\mathcal{L}$ is a trajectory $q$ such that $D \mathcal{L} (q)(w)=0$ for any \textit{variations} $w$ (i.e. $w \in \mathcal{C}^2([a,b],\Rr^d)$, $w(a)=w(b)=0$), where $D \mathcal{L} (q)(w)$ is the differential of $\mathcal{L}$ in $q$ along the direction $w$. Extremals of a Lagrangian functional can be characterized as solution of a differential equation of order $2$:

\begin{theo}[Variational principle]\label{vp}
Let $\mathcal{L}$ be a Lagrangian functional associated to the Lagrangian $L$ and let $q \in \mathcal{C}^{2} ([a,b],\Rr^d)$. Then, $q$ is an extremal of $\mathcal{L}$ if and only if $q$ is solution of the Euler-Lagrange equation given by: 
\begin{equation}\tag{$EL$}\label{el}
\forall t \in ]a,b[, \; \; \dis \frac{\partial L}{\partial x} (q(t),\dot{q}(t),t) - \frac{d}{dt} \left( \frac{\partial L}{\partial v} (q(t),\dot{q}(t),t) \right) = 0 . 
\end{equation}
\end{theo}

\noindent We now apply definitions of Section \ref{section1} on Lagrangian systems. 

\begin{prrop}\label{dde}
Let $L$ be a Lagrangian and let \eqref{el} be its associated Euler-Lagrange equation. The direct discrete embedding of \eqref{el} with respect to $FDE\sigma$ is given by: 
\begin{equation}\label{dd}
\frac{\partial L}{\partial x} (Q, -\sigma \Delta_{\sigma} Q, \tau) + \sigma \Delta_{\sigma} \left( \frac{\partial L}{\partial v} (Q,-\sigma \Delta_{\sigma} Q,\tau) \right) = 0, \; \; Q \in (\Rr^d)^{N+1} .
\end{equation}
\end{prrop}

\noindent We refer to \cite{bour} for a concrete example illustrating Theorem \ref{vp} and Proposition \ref{dde}. 

\section{Discrete embeddings and variational integrators of Lagrangian systems}\label{section2}
A direct discrete embedding is only based on the form of the differential operator which is dependent of the coordinates system and consequently is not intrinsic. Then, a natural question arises: \textit{what can be said about the conservation of intrinsic properties of a differential equation by a discrete embedding?} This paper is devoted to the conservation by discrete embeddings of the Lagrangian structure of Lagrangian systems. More precisely, Theorem \ref{vp} shows that \eqref{el} possesses a variational structure: \textit{the direct discrete embedding being a procedure mainly algebraic, does \eqref{dd} possess a variational structure too?} It is \textit{not always true}.

However, a variational integrator, developed in \cite[Chapter VI.6]{lubi} and in \cite{mars}, is a discretization of a Lagrangian system preserving its variational structure. Indeed, it is based on the discrete analogous of the variational principle on a discrete version of the associated Lagrangian functional. 

In our framework, the discretization of the Lagrangian functional is induced by giving a discrete embedding.

\subsection{Discrete Lagrangian functionals and discrete calculus of variations}\label{section21}
In this subsection, as an example, we are going to work exclusively in the framework of $FDE\sigma$. Giving $FDE\sigma$ induces the discretization of a Lagrangian functional as long as a quadrature formula is fixed in order to approximate integrals. We choose the usual $\sigma$-quadrature formula of Gauss: for a continuous function $f$ on $[a,b]$, we discretize $ \int_a^b f(t) dt$ by $h \sum_{k \in I_{\sigma}} f(t_k)$ where $I_+ = \{ 0,...,N-1 \} $ and $I_- = \{ 1,...,N \} $. \\

\noindent This process defines the \textit{Gauss finite differences embedding} denoted by $Gauss$-$FDE\sigma$. Such a choice allows to keep at the discrete level the following fundamental result:
\begin{equation*}
  \xymatrix@R=2cm@C=2.5cm {
  \dis \int _a ^b \dot{q}(t) dt = q(b) - q(a) \qquad \ar[r]^{Gauss\text{-}FDE\sigma} & \qquad h \dis \sum_{k \in I_{\sigma}} (-\sigma \Delta_{\sigma} Q )_k = Q_N - Q_0 .
  }
\end{equation*}

\begin{prrop}\label{dlf}
Let $\mathcal{L}$ be a Lagrangian functional associated to a Lagrangian $L$. The discrete Lagrangian functional associated to $\mathcal{L}$ with respect to $Gauss$-$FDE\sigma$ is given by:
$$ \fonction{\mathcal{L}^{\sigma}_h}{\left( \Rr^{d} \right)^{N+1}}{\Rr}{Q = (Q_k)_{k=0,...,N}}{h \displaystyle \sum_{k \in I_{\sigma}} L(Q_k,(-\sigma \Delta _{\sigma} Q)_k,t_k).} $$
\end{prrop}

\noindent Once the discrete version of the Lagrangian functional is formulated, we can develop a discrete calculus of variations on it: this leads to a variational integrator. Let $\mathcal{L}$ be a Lagrangian functional and $\mathcal{L}^{\sigma}_h$ the discrete Lagrangian functional associated with respect to $Gauss$-$FDE\sigma$. A \textit{discrete extremal} (or \textit{discrete critical point}) of $\mathcal{L}^{\sigma}_h$ is an element $Q$ in $(\Rr^{d})^{N+1}$ such that $D \mathcal{L}^{\sigma}_h (Q)(W)=0$ for any \textit{discrete variations} $W$ (i.e. $W \in (\Rr^{d})^{N+1}$, $W_0 = W_N =0$). Discrete extremals of $\mathcal{L}^{\sigma}_h$ can be characterized as solution of a system of equations:

\begin{theo}[Discrete variational principle]\label{dvp}
Let $\mathcal{L}^{\sigma}_h$ be the discrete Lagrangian functional associated to the Lagrangian $L$ with respect to $Gauss$-$FDE\sigma$. Then, $Q$ in $(\Rr^d)^{N+1}$ is a discrete extremal of $\mathcal{L}^{\sigma}_h$ if and only if $Q$ is solution of the following system of equations (called discrete Euler-Lagrange equation) given by:
\begin{equation}\label{elh}\tag{$EL^{\sigma}_h$}
\frac{\partial L}{\partial x} (Q, -\sigma \Delta_{\sigma} Q, \tau) - \sigma \Delta_{-\sigma} \left( \frac{\partial L}{\partial v} (Q,-\sigma \Delta_{\sigma} Q,\tau) \right) = 0, \; \; Q \in (\Rr^d)^{N+1} .
\end{equation}
\end{theo}

\noindent \eqref{elh} is obtained from \eqref{el} by variational integrator. Its variational origin allows us to say that it possesses a Lagrangian structure. Then, we have conservation at the discrete level of the Lagrangian structure by variational integrator. \\

\noindent Let us note that an asymmetry appears in \eqref{elh}: indeed, we have a composition between the two discrete operators $\Delta_+$ and $\Delta_-$. We notice that this asymmetry does not appear in the continuous space in \eqref{el}. 


\subsection{Problem of coherence of a discrete embedding}\label{section23}
Hence, defining a discrete embedding leads to two discrete versions of an Euler-Lagrange equation: the first one \eqref{dd} obtained by direct discrete embedding and the second one \eqref{elh} corresponding to a variational integrator. The direct discrete embedding is an algebraic procedure (respecting for example the law of semi-group of the differential operator $d/dt$). On the contrary, a variational integrator is mainly based on a dynamical approach via the extremals of a functional.

However, we are interested in the conservation at the discrete level of the Lagrangian structure of Lagrangian systems. We then propose the following definition: a discrete embedding is said to be \textit{coherent} if the two numerical schemes coincide. Precisely, a discrete embedding is coherent if it makes the following diagram commutative:

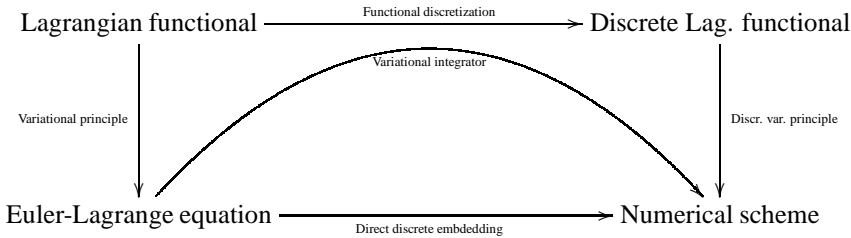
\begin{equation*}
  \xymatrix@R=2cm@C=4cm {
  \text{Lagrangian functional} \ar[r]^{\text{\begin{tiny} Functional discretization \end{tiny}}} \ar[d]_{\text{\begin{tiny} Variational principle \end{tiny}}} & \text{Discrete Lag. functional} \ar[d]^{\text{\begin{tiny} Discr. var. principle \end{tiny}}} \\
  \text{Euler-Lagrange equation} \ar[r]_{\text{\begin{tiny} Direct discrete embdedding \end{tiny}}} \ar@/^2.2cm/[r]_{\text{\begin{tiny} Variational integrator \end{tiny}}} & \text{Numerical scheme}
  }
\end{equation*}
~\\
\noindent Thus, a coherent discrete embedding provides a direct discrete version of a Lagrangian system preserving its Lagrangian structure. \\

\noindent The previous study leads to a default of coherence of $Gauss$-$FDE\sigma$. Indeed, algorithms obtained by direct discrete embedding \eqref{dd} and obtained by discrete variational principle \eqref{elh} do not coincide. The problem is to understand \textit{why there is not asymmetry appearing in the direct discrete embedding?} It seems that we miss dynamical informations in the formulation of Lagrangian systems at the continuous level which are pointed up in the discrete space with the asymmetric discrete operators $\left( -\sigma \Delta_{\sigma} \right)_{\sigma = \pm}$. Nevertheless, this default of coherence can be corrected using a different writing of the initial Euler-Lagrange equation. 

\subsection{Rewriting of the Euler-Lagrange equation and discrete embeddings}\label{section24}
The usual way to derive differential equations in Physics is to built a continuous model using discrete data. However, this process gives only an information in one direction of time. As a consequence, a discrete evaluation of the velocity corresponds in general at the continuous level to the evaluation of the right or left derivative. In general, we replace the right (or left) derivative by the classical derivative $d/dt$. However, this procedure assumes that the underlying solution is differentiable. This assumption is not only related to the regularity of the solutions but also to the reversibility of the systems (the right and left derivatives are equal). In this section, we introduce asymmetric Lagrangian systems which are obtained with functionals depending only on left or only on right derivatives. We prove in this case that $Gauss$-$FDE\sigma$ is coherent. 

\begin{defi}
For $f : [a,b] \longrightarrow \Rr ^d $ smooth enough function, we denote:
$$ \forall t \in ]a,b], \; \; d_- f(t) = \lim\limits_{h \to 0^+} \frac{f(t) - f(t-h)}{h}$$
and
$$ \forall t \in [a,b[, \; \; d_+ f(t) = \lim\limits_{h \to 0^+} \frac{f(t) - f(t+h)}{h} .$$
\end{defi}
\noindent Although we have $d_- f = -d_+ f = \dot{f}$ for a differentiable function $f$, it is interesting to use these notations in order to keep dynamical informations.

\begin{defi}
An asymmetric Lagrangian functional is an application:
$$ \fonction{\mathcal{L}^{\sigma}}{\mathcal{C}^{2}([a,b],\Rr^d)}{\Rr}{q}{\dis \int_{a}^{b} L(q(t),-\sigma d_{\sigma} q(t),t)dt} $$
where $L$ is a Lagrangian.
\end{defi}

\noindent Then, by calculus of variations, we obtain the following characterization of the extremals of an asymmetric Lagrangian functional:

\begin{theo}[Variational principle]\label{avp}
Let $\mathcal{L}^{\sigma}$ be an asymmetric Lagrangian functional associated to the Lagrangian $L$ and let $q \in \mathcal{C}^{2}([a,b],\Rr^d)$. Then, $q$ is an extremal of $\mathcal{L}^{\sigma}$ if and only if $q$ is solution of the asymmetric Euler-Lagrange equation:
\begin{equation}\label{ela}\tag{$EL^{\sigma}$}
\begin{array}{c}
\forall t \in ]a,b[, \; \; \dis \frac{\partial L}{\partial x} (q(t),-\sigma d_{\sigma} q(t),t) - \sigma d_{- \sigma} \left( \frac{\partial L}{\partial v} (q(t),-\sigma d_{\sigma} q(t),t) \right) = 0 . 
\end{array}
\end{equation}
\end{theo}

\noindent Hence, \eqref{ela} possesses a variational structure. \textit{Is it conserved by discrete embeddings?} In order to embed \eqref{ela}, we have to discretize two differential operators at the same time. We then define the following asymmetric version of $Gauss$-$FDE\sigma$:

\begin{defi}
We call the asymmetric version of $Gauss$-$FDE\sigma$ the definition of the following elements: the application $\mbox{disc}$, the $\sigma$-quadrature formula of Gauss and the discrete operators $\Delta_-$ and $\Delta_+$ discretizing respectively the operators $d_{-} $ and $d_{+}$.
\end{defi}

\begin{prrop}
The asymmetric version of $Gauss$-$FDE\sigma$ is a coherent discrete embedding. Indeed, the direct discrete embedding and the variational integrator of \eqref{ela} in the framework of the asymmetric $Gauss$-$FDE\sigma$ lead to the same numerical scheme: \eqref{elh}.
\end{prrop}

\noindent We notice that the rewriting \eqref{ela} of \eqref{el} provides additional dynamical informations which allows the asymmetric $Gauss$-$FDE\sigma$ to unify the algebraic and the dynamical approaches in the discretization of a Lagrangian system. Moreover, this rewriting can be justified by the fractional calculus as we will see in Section \ref{section3}.

\section{Discrete embeddings and variational integrators of fractional Lagrangian systems}\label{section3}
\subsection{Fractional derivatives and fractional Lagrangian systems}\label{section31}
Fractional calculus is the generalization of the derivative notion to real orders. We refer to \cite{hilf}, \cite{samk} for many different ways generalizing this notion. For the whole paper, we fix $0 < \alpha < 1$ and for any $r \in \Nn^{\ast}$, we denote by $\alpha_r = (-\alpha)(1-\alpha)...(r-1-\alpha)/r!$ and $\alpha_0 = 1$. We are going to use the classical notions of Gr\"unwald-Letnikov. The following definition is extracted from \cite{podl}.

\begin{defi} 
Let $f$ be an element of $\mathcal{C}^{1} ([a,b],\Rr^d)$. The Gr\"unwald-Letnikov fractional left derivative of order $\alpha$ with inferior limit $a$ of $f$ is:
$$\forall t \in ]a,b], \; \; D^{\alpha}_{-} f(t) = \lim\limits_{\substack{h \to 0 \\ nh = t-a}} \dfrac{1}{h^{\alpha}} \displaystyle \sum_{r=0}^{n} \alpha_r f(t-rh) $$ 
and the Gr\"unwald-Letnikov fractional right derivative of order $\alpha$ with superior limit $b$ of $f$ is:
$$\forall t \in [a,b[, \; \; D^{\alpha}_{+} f(t) = \lim\limits_{\substack{h \to 0 \\ nh = b-t}} \dfrac{1}{h^{\alpha}} \displaystyle \sum_{r=0}^{n} \alpha_r f(t+rh) . $$ 
\end{defi}

\noindent Recently, an important activity has been devoted to fractional Lagrangian systems for the purpose of optimal control, mechanics, engineering and Physics, \cite{agra}, \cite{agra3}, \cite{fred}. We recall definitions and results concerning these fractional systems, we refer to \cite{agra} for a detailed study and for a detailed proof of Theorem \ref{vpf}. 

\begin{defi}
A fractional Lagrangian functional of order $\alpha$ is an application defined by:
\begin{equation*}
\fonction{\mathcal{L}^{\sigma,\alpha}}{\mathcal{C}^2 ([a,b],\Rr ^d)}{\Rr}{q}{\dis \int_{a}^{b} L(q(t),-\sigma D^{\alpha}_{\sigma} q(t),t)dt}
\end{equation*}
where $L$ is a Lagrangian.
\end{defi}
\noindent We can give a characterization of extremals of a fractional Lagrangian functional as solutions of a fractional differential equation:

\begin{theo}[Variational principle]\label{vpf}
Let $\mathcal{L}^{\sigma,\alpha}$ be a fractional Lagrangian functional of order $\alpha$ associated to the Lagrangian $L$ and let $q$ be an element of $\mathcal{C}^{2}([a,b],\Rr^d)$. Then, $q$ is an extremal of $\mathcal{L}^{\sigma,\alpha}$ if and only if $q$ is solution of the fractional Euler-Lagrange equation:
\begin{equation}\label{elf}\tag{$EL^{\sigma,\alpha}$}
\forall t \in ]a,b[, \; \; \frac{\partial L}{\partial x} (q(t),-\sigma D^{\alpha}_{\sigma} q(t),t) - \sigma D^{\alpha}_{-\sigma} \left( \frac{\partial L}{\partial v} (q(t),-\sigma D^{\alpha}_{\sigma} q(t),t) \right) = 0 . 
\end{equation}
\end{theo}

\noindent We refer to \cite{agra} for a detailed proof. Hence, in the fractional case, we find an asymmetry again making a link with the asymmetric rewriting of \eqref{el} into \eqref{ela}. \\

\noindent As in the classical case, we conclude that \eqref{elf} possesses a Lagrangian structure and we are iterested by its conservation at the discrete level by discrete embeddings.

\subsection{Discrete embeddings of fractional Lagrangian systems}\label{section32}
There exist many studies concerning the discretization of fractional differential equations but without the point of view of discrete embeddings. We refer to \cite{agra3}, \cite{bast}. By referring to the notion of Gr\"unwald-Letnikov \cite{dubo}, we give the following definition: 

\begin{defi}
The Gauss Gr\"unwald-Letnikov embedding denoted by $Gauss$-$GLE\sigma$ is the definition of the following elements: the application $\mbox{disc}$, the $\sigma$-quadrature formula of Gauss and the discrete operators
$$ \fonction{\Delta^{\alpha}_{-}}{(\Rr ^d )^{N+1}}{(\Rr ^d )^{N}}{Q = (Q_k)_{k=0,...,N}}{\left( \dfrac{1}{h^{\alpha}} \displaystyle \sum_{r=0}^{k} \alpha_r Q_{k-r} \right)_{k=1,..,N} } $$
and
$$ \fonction{\Delta^{\alpha}_{+}}{(\Rr ^d )^{N+1}}{(\Rr ^d )^{N}}{Q = (Q_k)_{k=0,...,N}}{\left( \dfrac{1}{h^{\alpha}} \displaystyle \sum_{r=0}^{N-k} \alpha_r Q_{k+r} \right)_{k=0,..,N-1} .} $$
These discrete operators are respectively the discrete versions of $D^{\alpha}_{-} $ and $D^{\alpha}_{+}$. 
\end{defi}

\noindent We are first interested in the variational integrator of \eqref{elf} in the framework of $Gauss$-$GLE\sigma$. Giving $Gauss$-$GLE\sigma$ allows us to formulate the discrete version of a fractional Lagrangian functional:

\begin{prrop}
Let $\mathcal{L}^{\sigma,\alpha}$ be the fractional Lagrangian functional associated to the Lagrangian $L$. The discrete fractional Lagrangian functional associated to $\mathcal{L}^{\sigma,\alpha}$ with respect to $Gauss$-$GLE\sigma$ is given by:
$$ \fonction{\mathcal{L}^{\sigma,\alpha}_h}{\left( \Rr^{d} \right)^{N+1}}{\Rr}{Q = (Q_k)_{k=0,...,N}}{h \dis \sum_{k \in I_{\sigma}} L(Q_k,(-\sigma \Delta^{\alpha}_{\sigma} Q)_k,t_k).} $$
\end{prrop}
\noindent Then, discrete extremals of the discrete fractional Lagrangian functional can be characterized as solutions of a system of equations:

\begin{theo}[Discrete variational principle]\label{dflap}
Let $\mathcal{L}^{\sigma,\alpha}_h$ be a discrete fractional Lagrangian functional associated to the Lagrangian $L$ with respect to $Gauss$-$GLE\sigma$. Then, $Q$ in $(\Rr^d)^{N+1}$ is a discrete extremal of $\mathcal{L}^{\sigma,\alpha}_h$ if and only if $Q$ is solution of the following system of equations, called the discrete fractional Euler-Lagrange equation:
\begin{equation}\label{elfh}\tag{$EL^{\sigma,\alpha}_h$}
\frac{\partial L}{\partial x} (Q,-\sigma \Delta ^{\alpha}_{\sigma} Q,\tau) - \sigma \Delta ^{\alpha}_{-\sigma} ( \frac{\partial L}{\partial v} (Q,-\sigma \Delta ^{\alpha}_{\sigma} Q,\tau)) = 0, \; \; Q \in (\Rr^d)^{N+1}.
\end{equation}
\end{theo}

\noindent We conclude with the following proposition:

\begin{prrop}
$Gauss$-$GLE\sigma$ is a coherent discrete embedding. Indeed, the direct discrete embedding and the variational integrator of \eqref{elf} in the framework of $Gauss$-$GLE\sigma$ lead to the same numerical scheme: \eqref{elfh}.
\end{prrop}


%
  
\bibliographystyle{acmurl}
\bibliography{proceedingbib}

\begin{address}
  LMA, University of Pau -  Postal address IPRA BP 1155 Pau Cedex (France) \\
   \texttt{bourdin.l@etud.univ-pau.fr}
 \end{address}

\end{document}